# FROM CONTINUOUS-TIME DESIGN TO SAMPLED-DATA DESIGN OF NONLINEAR OBSERVERS

Iasson Karafyllis[*] and Costas Kravaris[**]

[*]Department of Environmental Engineering, Technical University of Crete,
73100, Chania, Greece
email: ikarafyl@enveng.tuc.gr

[**]Department of Chemical Engineering, University of Patras
26500 Patras, Greece
email: kravaris@chemeng.upatras.gr

**Abstract**

In this work, a sampled-data nonlinear observer is designed using a continuous-time design coupled with an inter-sample output predictor. The proposed sampled-data observer is a hybrid system. It is shown that under certain conditions, the robustness properties of the continuous-time design are inherited by the sampled-data design, as long as the sampling period is not too large. The approach is applied to linear systems and to triangular globally Lipschitz systems.

## 1. Introduction

The problem of designing sampled-data nonlinear observers is a very challenging problem that has attracted a lot of attention in the literature. Continuous-time nonlinear observer designs [12,18,20,25,26,27] are meant to be used only for very small sampling periods, whereas their potential "redesign" for the purpose of digital implementation, even though straightforward and popular for linear systems [5], poses significant challenges in the nonlinear case. For this reason, the main line of attack has been through the use of an exact or approximate discrete-time description of the dynamics as the starting point for observer design [3,4,7,8,9,10,11,14,19,21,22,26,28]. This is a reasonable point of view, but faces two important difficulties:
(i) from the moment that the continuous-time system description is abandoned and is substituted by a discrete-time description, the inter-sample dynamic behavior is lost
(ii) any errors in the sampling schedule, get transferred into errors in the discrete-time description
As a consequence, available design methods (i) do not provide an explicit estimate of the error in between two consecutive sampling times and (ii) do not account for perturbations of the sampling schedule. Moreover, due to observability issues, the magnitude of the sampling period cannot be arbitrary (see [1,26]).

Finally, optimization-based approaches for nonlinear observer design [2,13,23,24,29] are also based on a discrete-time description of the dynamics and therefore share all the above difficulties, but, because they utilize a large number of measurements, offer the advantage of reduced sensitivity to measurement errors at the expense of higher memory requirements and computational cost .

A hybrid observer design approach was recently proposed in [6], which bears similarities to the above-mentioned optimization-based approach, but the hybrid nature of their observer offers certain advantages. In the present work, our proposed sampled-data observer will also be a hybrid system; however, it will directly emerge from a continuous-time design of a nonlinear observer.

Consider a single-output continuous-time system:

$$\dot{x} = f(x), x \in \Re^n$$
$$y = h(x), y \in \Re \tag{1.1}$$

$f \in C^1(\Re^n; \Re^n)$, $h \in C^2(\Re^n; \Re)$ with $f(0) = 0$, $h(0) = 0$. For this system, suppose that a continuous-time observer design is available

$$\dot{z} = F(z, y), z \in \Re^k$$
$$\hat{x} = \Psi(z), \hat{x} \in \Re^n \qquad (1.2)$$

where $F \in C^1(\Re^k \times \Re; \Re^k)$, $\Psi \in C^1(\Re^k; \Re^n)$ with $F(0,0) = 0$, $\Psi(0) = 0$.

The question is whether this design would still be useful in the presence of sampled measurements $y(ih)$, $i=0,1,\ldots$, where $h$ is the sampling period, or more generally, at some countable set of time instants $\pi = \{\tau_i\}_{i=0}^{\infty}$, not necessarily uniformly spaced, but satisfying $0 < \tau_{i+1} - \tau_i \le r$ for all $i = 0,1,\ldots$ for some $r > 0$.

The present work has been motivated by the intuitive expectation that a continuous-time nonlinear observer design would still be useful in the presence of "medium-size" sampling periods, as long as special care is taken in the time-interval between measurements. Holding the most recent measurement (zero-order hold) is not the most intuitively meaningful strategy; instead, the model (1.1) could be used to predict the evolution of the output, up until the new measurement is received. In particular, in the present paper, we propose a sampled-data observer consisting of the continuous-time observer, coupled with an output predictor for the time interval between two consecutive measurements:

$$\dot{z}(t) = F(z(t), w(t)), t \in [\tau_i, \tau_{i+1})$$
$$\dot{w}(t) = L_f h(\Psi(z(t))), t \in [\tau_i, \tau_{i+1})$$
$$w(\tau_{i+1}) = y(\tau_{i+1}) \qquad (1.3)$$
$$(z(t), w(t)) \in \Re^k \times \Re$$
$$\hat{x}(t) = \Psi(z(t)), \hat{x} \in \Re^n$$

Figure 1 depicts the structure of the sampled-data observer (1.3) compared to the continuous-time observer (1.2). The sampled-data observer uses the continuous-time observer as a key ingredient, coupled with an inter-sample output predictor. The latter is initialized by the most recent measurement and integrates the rate of change of the output calculated by the model ($L_f h(x) := \nabla h(x) f(x)$ is the Lie derivative of the output map).

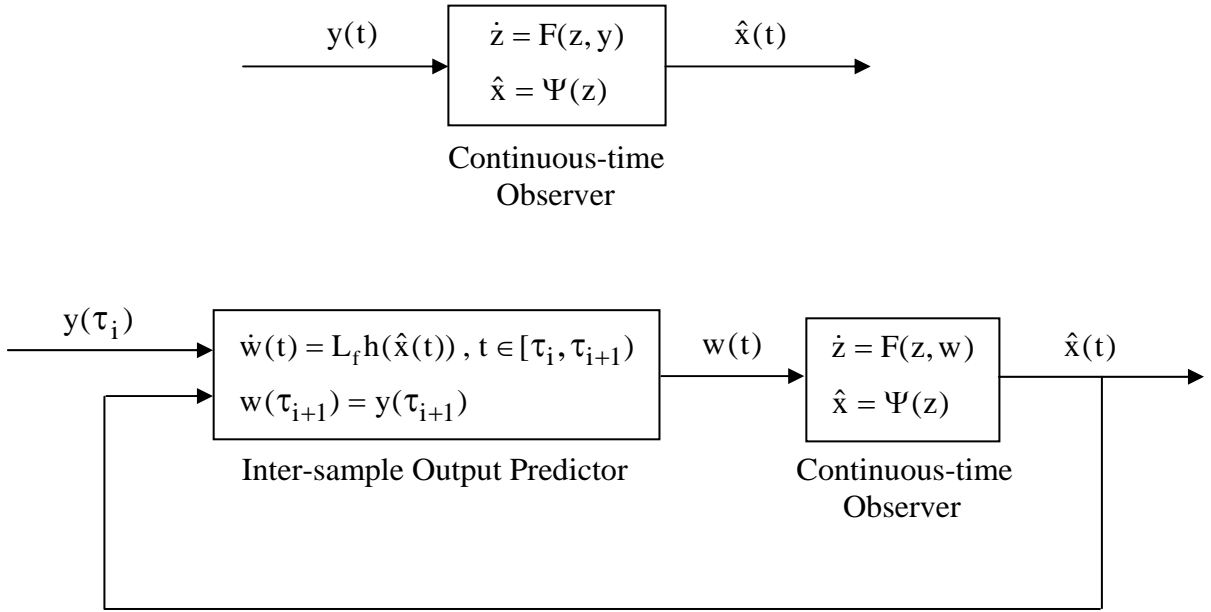

**Figure 1:** Continuous-time observer (1.2) (top) versus sampled-data observer (1.3) (bottom).

It is important point out that the entire system (1.1) with (1.3) is a hybrid system, which does not satisfy the classical semigroup property. However, the weak semigroup property holds (see [15,16]) and consequently it can be analyzed using the recent results in [15,16,17].

The main result of the present paper is that the properties of the observer (1.2) under continuous measurement are inherited by the observer (1.3) under sampled measurements, as long as the sampling period is not too large.

The paper is organized as follows. First, basic notations are defined, followed by definitions of notions of robust observer for forward complete dynamic systems, which will play a central role in the rest of the paper. Next, the main result of the paper is presented and proven. Finally, the main result is applied to two classes of systems (i) triangular globally Lipschitz systems and (ii) linear systems, leading to more concrete sampled-data observer designs in these classes of systems.

**Notations** Throughout this paper we adopt the following notations:

* Let $A \subseteq \Re^n$ be an open set. By $C^0(A;\Omega)$, we denote the class of continuous functions on $A$, which take values in $\Omega \subseteq \Re^k$. By $C^l(A;\Omega)$, where $l \in \{1,2,...\}$, we denote the class of continuous functions on $A$ with continuous derivatives of order $l$, which take values in $\Omega$.

* For a vector $x \in \Re^n$ we denote by $|x|$ its usual Euclidean norm and by $x'$ its transpose. By $|A| := \sup\{|Ax|; x \in \Re^n, |x| = 1\}$ we denote the induced norm of a matrix $A \in \Re^{m \times n}$ and $I$ denotes the identity matrix. By $B = diag(b_1,...,b_n)$ we denote the diagonal matrix $B \in \Re^{n \times n}$ with $b_1,...,b_n$ in its diagonal.

* By $\Re^+$ we denote the set of non-negative real numbers.

* We denote by $K^+$ the class of positive $C^0$ functions defined on $\Re^+$. We say that a non-decreasing continuous function $\gamma : \Re^+ \to \Re^+$ is of class $N$ if $\gamma(0) = 0$. We say that a function $\rho : \Re^+ \to \Re^+$ is positive definite if $\rho(0) = 0$ and $\rho(s) > 0$ for all $s > 0$. By $K$ we denote the set of positive definite, increasing and continuous functions. We say that a positive definite, increasing and continuous function $\rho : \Re^+ \to \Re^+$ is of class $K_\infty$ if $\lim_{s \to +\infty} \rho(s) = +\infty$. By $KL$ we denote the set of all continuous functions $\sigma = \sigma(s,t) : \Re^+ \times \Re^+ \to \Re^+$ with the properties: (i) for each $t \geq 0$ the mapping $\sigma(\cdot,t)$ is of class $K$; (ii) for each $s \geq 0$, the mapping $\sigma(s,\cdot)$ is non-increasing with $\lim_{t \to +\infty} \sigma(s,t) = 0$.

* Let $D \subseteq \Re^l$ be a non-empty set and $I \subseteq \Re^+$ an interval. By $L^\infty_{loc}(I;D)$ we denote the class of all Lebesgue measurable and locally bounded mappings $d : \Re^+ \to D$. Notice that by $\sup_{\tau \in [0,t]} |d(\tau)|$ we do not mean the essential supremum of $d : \Re^+ \to D$ on $[0,t]$ but the actual supremum of $d : \Re^+ \to D$ on $[0,t]$.

* Let $f \in C^1(\Re^n;\Re^n)$, $h \in C^1(\Re^n;\Re)$. By $L_f h(x) := \nabla h(x) f(x)$ we denote the Lie derivative of the function $h \in C^1(\Re^n;\Re)$ along the vector field $f \in C^1(\Re^n;\Re^n)$.

## 2. Basic Notions

In the present work we study systems of the form (1.1) under the following hypotheses:

**(H)** System (1.1) is Robustly Forward Complete (see [15]), i.e., there exist functions $\mu \in K^+$ and $a \in K_\infty$ such that for every $x_0 \in \Re^n$ the solution $x(t)$ of (1.1) with initial condition $x(0) = x_0$ satisfies

$$|x(t)| \leq \mu(t) a(|x_0|), \quad \forall t \geq 0 \tag{2.1}$$

The following definition of the notion of robust observer for system (1.1) with respect to measurement errors is crucial to the development of the main results of the present work.

**Definition 2.1:** *Consider the following system*

$$\dot{z} = F(z,y), z \in \Re^k \\ \hat{x} = \Psi(z), \hat{x} \in \Re^n \tag{2.2}$$

where $F \in C^1(\Re^k \times \Re; \Re^k)$, $\Psi \in C^1(\Re^k; \Re^n)$ with $F(0,0) = 0$, $\Psi(0) = 0$. System (2.2) is called a **robust observer** for system (1.1) with respect to measurement errors, if the following conditions are met:

i) there exist functions $\sigma \in KL$, $\gamma, p \in N$, $\mu \in K^+$ and $a \in K_\infty$ such that for every $(x_0, z_0) \in \Re^n \times \Re^k$ and $v \in L_{loc}^\infty(\Re^+; \Re)$, the solution $(x(t), z(t))$ of

$$\dot{x} = f(x)$$
$$\dot{z} = F(z, h(x) + v) \quad (2.3)$$
$$\hat{x} = \Psi(z)$$

with initial condition $(x(0), z(0)) = (x_0, z_0)$ corresponding to $v \in L_{loc}^\infty(\Re^+; \Re)$ satisfies the following estimates:

$$|\hat{x}(t) - x(t)| \leq \sigma(|x_0| + |z_0|, t) + \sup_{0 \leq \tau \leq t} \gamma(|v(\tau)|), \quad \forall t \geq 0 \quad (2.4a)$$

$$|z(t)| \leq \mu(t)\left[a(|x_0| + |z_0|) + \sup_{0 \leq \tau \leq t} p(|v(\tau)|)\right], \quad \forall t \geq 0 \quad (2.4b)$$

ii) for every $x_0 \in \Re^n$ there exists $z_0 \in \Re^k$ such that the solution $(x(t), z(t))$ of (2.3) with initial condition $(x(0), z(0)) = (x_0, z_0)$ corresponding to $v \equiv 0$, satisfies $x(t) = \Psi(z(t))$ for all $t \geq 0$.

**Remark 2.2:** If system (2.2) is a robust observer for system (1.1) with respect to measurement errors, then system (2.3) with output $Y = \Psi(z) - x$ satisfies the Uniform Input-to-Output Stability property from the input $v \in L_{loc}^\infty(\Re^+; \Re)$ with gain $\gamma \in N$ (see [17]).

We next define the corresponding notion of robust sampled-data observer. Notice that contrary to usual observers for which the output signal $y(t)$ of system (1.1) is available on-line, a sampled-data observer uses only the output values $y(\tau_i)$ at certain time instances $\pi = \{\tau_i\}_{i=0}^\infty$ with $0 < \tau_{i+1} - \tau_i \leq r$ for all $i = 0,1,...$. The number $r > 0$ is called the *upper diameter* of the sampling partition.

**Definition 2.3:** *The system*

$$\dot{z}(t) = g(z(t), z(\tau_i), y(\tau_i)), \, t \in [\tau_i, \tau_{i+1})$$
$$z(\tau_{i+1}) = G\left(\lim_{t \to \tau_{i+1}^-} z(t), y(\tau_{i+1})\right) \quad (2.5)$$
$$\hat{x}(t) = \Psi(z(t))$$

where $g \in C^1(\Re^k \times \Re^k \times \Re; \Re^k)$, $G \in C^0(\Re^k \times \Re; \Re^k)$, $\Psi \in C^1(\Re^k; \Re^n)$ with $g(0,0,0) = 0$, $G(0,0) = 0$, $\Psi(0) = 0$, *is called a **robust sampled-data observer** for (1.1) with respect to measurement errors, if the following conditions are met:*

i) there exist functions $\sigma \in KL$, $\gamma, p \in N$, $\mu \in K^+$ and $a \in K_\infty$ such that for every $(x_0, z_0, d, v) \in \Re^n \times \Re^k \times L_{loc}^\infty(\Re^+; \Re^+) \times L_{loc}^\infty(\Re^+; \Re)$, the solution $(x(t), z(t))$ of

$$\dot{x}(t) = f(x(t))$$
$$\dot{z}(t) = g(z(t), z(\tau_i), y(\tau_i) + v(\tau_i)), \, t \in [\tau_i, \tau_{i+1})$$
$$z(\tau_{i+1}) = G\left(\lim_{t \to \tau_{i+1}^-} z(t), y(\tau_{i+1}) + v(\tau_{i+1})\right) \quad (2.6)$$
$$\tau_{i+1} = \tau_i + r \exp(-d(\tau_i))$$
$$\hat{x}(t) = \Psi(z(t))$$

with initial condition $(x(0), z(0)) = (x_0, z_0)$ corresponding to $d \in L_{loc}^\infty(\Re^+; \Re^+)$, $v \in L_{loc}^\infty(\Re^+; \Re)$ satisfies the following estimates:

$$|\hat{x}(t) - x(t)| \leq \sigma(|x_0| + |z_0|, t) + \sup_{0 \leq \tau \leq t} \gamma(|v(\tau)|), \quad \forall t \geq 0 \qquad (2.7a)$$

$$|z(t)| \leq \mu(t)\left[a(|x_0| + |z_0|) + \sup_{0 \leq \tau \leq t} p(|v(\tau)|)\right], \quad \forall t \geq 0 \qquad (2.7b)$$

ii) for every $x_0 \in \Re^n$ there exists $z_0 \in \Re^k$ such that for all $d \in L_{loc}^\infty(\Re^+; \Re^+)$ the solution $(x(t), z(t))$ of (2.6) with initial condition $(x(0), z(0)) = (x_0, z_0)$ corresponding to $d \in L_{loc}^\infty(\Re^+; \Re^+)$ and $v \equiv 0$, satisfies $x(t) = \Psi(z(t))$ for all $t \geq 0$.

**Remark 2.4:** For each $(t_0, x_0, z_0) \in \Re^+ \times \Re^n \times \Re^k$ and for each $d \in L_{loc}^\infty(\Re^+; \Re^+)$, $v \in L_{loc}^\infty(\Re^+; \Re)$ the solution $(x(t), z(t))$ of (2.6) with initial condition $(x(0), z(0)) = (x_0, z_0)$ corresponding to $d \in L_{loc}^\infty(\Re^+; \Re^+)$, $v \in L_{loc}^\infty(\Re^+; \Re)$ is produced by the following algorithm:

Step $i$:
1) Given $\tau_i$ and $d \in L_{loc}^\infty(\Re^+; \Re^+)$, calculate $\tau_{i+1}$ using the equation $\tau_{i+1} = \tau_i + r \exp(-d(\tau_i))$,
2) Compute the state trajectory $(x(t), z(t))$, $t \in [\tau_i, \tau_{i+1}]$ as the solution of the differential equation $\dot{x}(t) = f(x(t))$ and $\dot{z}(t) = g(z(t), z(\tau_i), h(x(\tau_i)) + v(\tau_i))$,
3) Calculate $z(\tau_{i+1})$ using the equation $z(\tau_{i+1}) = G\left(\lim_{t \to \tau_{i+1}^-} z(t), h(x(\tau_{i+1})) + v(\tau_{i+1})\right)$,
4) Compute the output trajectory $Y(t)$, $t \in [\tau_i, \tau_{i+1}]$ using the equation $Y(t) = \Phi(z(t)) - x(t)$

For $i = 0$ we take $\tau_0 = t_0$ and $x(\tau_0) = x_0$ (initial condition). Hybrid systems of the form (2.6) were studied in [15,16,17], where the weak semigroup property for such systems was exploited. Taking into account hypothesis (H) for system (1.1), regularity properties of the right hand-sides of (2.5) and using the results of [15], we may conclude that

(i) system (2.6) has the **"Boundedness-Implies-Continuation" (BIC)** property, i.e., for each $(t_0, x_0, z_0, d, v) \in \Re^+ \times \Re^n \times \Re^k \times L_{loc}^\infty(\Re^+; \Re^+) \times L_{loc}^\infty(\Re^+; \Re)$, there exists $t_{max} \in (t_0, +\infty]$ (the maximal existence time of the solution) such that the solution $(x(t), z(t))$ of (2.6) with initial condition $(x(t_0), z(t_0)) = (x_0, z_0)$ corresponding to $d \in L_{loc}^\infty(\Re^+; \Re^+)$, $v \in L_{loc}^\infty(\Re^+; \Re)$ exists for all $t \in [t_0, t_{max})$. In addition, if $t_{max} < +\infty$ then for every $C > 0$ there exists $t \in [t_0, t_{max})$ with $|z(t)| > C$.

(ii) $0 \in \Re^n \times \Re^k$ is a **robust equilibrium point** from the input $(d, v) \in L_{loc}^\infty(\Re^+; \Re^+) \times L_{loc}^\infty(\Re^+; \Re)$, i.e., for every $\varepsilon > 0$, $T \in \Re^+$ there exists $\delta := \delta(\varepsilon, T) > 0$ such that for all $(t_0, x_0, z_0, d, v) \in \Re^+ \times \Re^n \times \Re^k \times L_{loc}^\infty(\Re^+; \Re^+) \times L_{loc}^\infty(\Re^+; \Re)$, with $|x_0| + |z_0| + \sup_{t \geq 0} d(t) + \sup_{t \geq 0} |v(t)| < \delta$ it holds that the solution $(x(t), z(t))$ of (2.6) with initial condition $(x(t_0), z(t_0)) = (x_0, z_0)$ corresponding to $(d, v) \in L_{loc}^\infty(\Re^+; \Re^+) \times L_{loc}^\infty(\Re^+; \Re)$ exists for all $t \in [t_0, t_0 + T]$ and

$$\sup\{|(x(t), z(t))|; \ t \in [t_0, t_0 + T], t_0 \in [0, T]\} < \varepsilon$$

(iii) system (2.6) is **autonomous**, i.e., for each $(t_0, x_0, z_0, d, v) \in \Re^+ \times \Re^n \times \Re^k \times L_{loc}^\infty(\Re^+; \Re^+) \times L_{loc}^\infty(\Re^+; \Re)$, $t \geq t_0$ and for each $\theta \in (-\infty, t_0]$ it holds that the solution $(x(t), z(t))$ of (2.6) with initial condition $(x(t_0), z(t_0)) = (x_0, z_0)$ corresponding to $d \in L_{loc}^\infty(\Re^+; \Re^+)$, $v \in L_{loc}^\infty(\Re^+; \Re)$ coincides with $(\tilde{x}(t - \theta), \tilde{z}(t - \theta))$ where $(\tilde{x}(t), \tilde{z}(t))$ is the solution of (2.6) with initial condition $(\tilde{x}(t_0 - \theta), \tilde{z}(t_0 - \theta)) = (x_0, z_0)$ corresponding to $P_\theta d \in L_{loc}^\infty(\Re^+; \Re^+)$ and $P_\theta v \in L_{loc}^\infty(\Re^+; \Re^+)$, where $(P_\theta d)(t) = d(t + \theta)$ and $(P_\theta v)(t) = v(t + \theta)$ for all $t + \theta \geq 0$.

**Remark 2.5:** If system (2.5) is a robust sampled-data observer for system (1.1) with respect to measurement errors, then system (2.6) with output $Y = \Psi(z) - x$ satisfies the Uniform Input-to-Output Stability property from the input $v \in L_{loc}^{\infty}(\Re^+; \Re)$ with gain $\gamma \in N$ (see [17], where the notion of Uniform Input-to-Output Stability (UIOS) was defined for hybrid systems ).

**Remark 2.6:** The reader should notice that the sampling period is allowed to be time-varying. The factor $\exp(-d(\tau_i)) \leq 1$, with $d(t) \geq 0$ some non-negative function of time, is an uncertainty factor in the end-point of the sampling interval. Proving stability for any non-negative input $d \in L_{loc}^{\infty}(\Re^+; \Re^+)$ will guarantee stability for all sampling schedules with $\tau_{i+1} - \tau_i \leq r$ (robustness to perturbations of the sampling schedule). To understand the importance of robustness to perturbations of the sampling schedule, consider the following situation. Suppose that hardware limitations restrict the sampling period to be $1s$. If we manage to design a sampled-data observer with $r \geq 2s$, then the application of the sampled-data observer will guarantee convergence of the state estimates even if we "miss measurements" or if we have "delayed measurements" (for example, due to improper operation of the sensor). In such a case robustness to perturbations of the sampling schedule becomes critical. The introduction of the factor $\exp(-d(\tau_i)) \leq 1$ is a mathematical way of introducing perturbations to the sampling schedule; however, it is not unique. Other ways of introducing perturbations of the sampling schedule can be considered.

## 3. Main Results

We are now in a position to state our main result.

**Theorem 3.1:** *Consider system (1.1) under hypothesis (H) and suppose that system (2.2) is a robust observer for system (1.1) with respect to measurement errors. Moreover, suppose that there exists a constant $K \geq 0$ and a function $\bar{\sigma} \in KL$ such that for every $(x_0, z_0) \in \Re^n \times \Re^k$ and $v \in L_{loc}^{\infty}(\Re^+; \Re)$, the solution $(x(t), z(t))$ of (2.3) satisfies the following estimate:*

$$\left| L_f h(\Psi(z(t))) - L_f h(x(t)) \right| \leq \bar{\sigma}\left( |x_0| + |z_0|, t \right) + K \sup_{0 \leq \tau \leq t} |v(\tau)|, \quad \forall t \geq 0 \tag{3.1}$$

*Finally, suppose that $rK < 1$, where $r > 0$ is the upper diameter of the sampling partition and $K \geq 0$ is the constant involved in estimate (3.1). Then (1.3) is a robust sampled-data observer for system (1.1) with respect to measurement errors.*

**Remark 3.2:** The reader should notice the structural differences between the continuous time observer (2.2) and the sampled-data observer (1.3), which are shown in Figures 1 and 2. The sampled-data observer uses the estimate of the state $\hat{x}(t)$ and the measurement $y(\tau_i)$ in order to generate an additional signal $w(t)$: the signal $w(t)$ will approximate the output signal $y(t)$ and actually replaces the output signal $y(t)$ in the observer.

**Proof of Theorem 3.1:** By virtue of Definition 2.3, it suffices to show that the following hybrid system:

$$\begin{aligned}
\dot{x}(t) &= f(x(t)) \\
\dot{z}(t) &= F(z(t), w(t)), \, t \in [\tau_i, \tau_{i+1}) \\
\dot{w}(t) &= L_f h(\Psi(z(t))), \, t \in [\tau_i, \tau_{i+1}) \\
w(\tau_{i+1}) &= y(\tau_{i+1}) + v(\tau_{i+1}) \\
\tau_{i+1} &= \tau_i + r \exp(-d(\tau_i)) \\
Y(t) &= \Psi(z(t)) - x(t)
\end{aligned} \tag{3.2}$$

satisfies the Uniform Input-to-Output Stability property from the input $v \in L_{loc}^{\infty}(\Re^+; \Re)$, i.e., it suffices to show that system (3.2) is Robustly Forward Complete from the input $v \in L_{loc}^{\infty}(\Re^+; \Re)$, $0 \in \Re^n \times \Re^k \times \Re$ is a robust equilibrium point from the input $v \in L_{loc}^{\infty}(\Re^+; \Re)$ (see [15,16,17]) and that there exist functions $\sigma \in KL$, $\mu \in K^+$ $\tilde{\gamma} \in N$, and $a \in K_\infty$ such that for every $(x_0, z_0, w_0, d, v) \in \Re^n \times \Re^k \times \Re \times L_{loc}^\infty(\Re^+; \Re^+) \times L_{loc}^\infty(\Re^+; \Re)$ the solution $(x(t), z(t), w(t))$ of (3.2) with initial condition $(x(0), z(0), w(0)) = (x_0, z_0, w_0)$ corresponding to $d \in L_{loc}^\infty(\Re^+; \Re^+)$, $v \in L_{loc}^\infty(\Re^+; \Re)$ satisfies:

$$|Y(t)| \leq \sigma(|x_0| + |z_0| + |w_0|, t) + \sup_{0 \leq \tau \leq t} \tilde{\gamma}(|v(\tau)|), \quad \forall t \geq 0 \tag{3.3}$$

The reader should notice that for every $x_0 \in \Re^n$ there exists $(z_0, w_0) \in \Re^k \times \Re$ with $w_0 = h(x_0)$ such that for all $d \in L_{loc}^\infty(\Re^+; \Re^+)$ the solution $(x(t), z(t), w(t))$ of (3.2) with initial condition $(x(0), z(0), w(0)) = (x_0, z_0, w_0)$ corresponding to $d \in L_{loc}^\infty(\Re^+; \Re^+)$, $v \equiv 0$, satisfies $x(t) = \Psi(z(t))$ for all $t \geq 0$.

Since system (2.2) is a robust observer for system (1.1) with respect to measurement errors and since hypothesis (H) holds, it follows from (2.1), (2.4a,b) and (3.1) that for every $(x_0, z_0, w_0, d) \in \Re^n \times \Re^k \times \Re \times L_{loc}^\infty(\Re^+; \Re^+)$ the solution $(x(t), z(t), w(t)) \in \Re^n \times \Re^k \times \Re$ of (3.2) with initial condition $(x(0), z(0), w(0)) = (x_0, z_0, w_0)$ satisfies the following estimates:

$$|Y(t)| \leq \sigma(|x_0| + |z_0|, t) + \sup_{0 \leq \tau \leq t} \gamma(|w(\tau) - h(x(\tau))|), \quad \forall t \in [0, t_{\max}) \tag{3.4}$$

$$|L_f h(\Psi(z(t))) - L_f h(x(t))| \leq \bar{\sigma}(|x_0| + |z_0|, t) + K \sup_{0 \leq \tau \leq t} |w(\tau) - h(x(\tau))|, \quad \forall t \in [0, t_{\max}) \tag{3.5}$$

$$|z(t)| + |x(t)| \leq \mu(t) \left[ a(|x_0| + |z_0|) + \sup_{0 \leq \tau \leq t} p(|w(\tau) - h(x(\tau))|) \right], \quad \forall t \in [0, t_{\max}) \tag{3.6}$$

for appropriate functions $\sigma, \bar{\sigma} \in KL$, $\gamma, p \in N$, $\mu \in K^+$ and $a \in K_\infty$, where $t_{\max} \in (0, +\infty]$ is the maximal existence time of the solution. Let $\pi = \{\tau_i\}_{i=0}^\infty$ be the partition of $\Re^+$ generated by the recursive formula $\tau_{i+1} = \tau_i + r \exp(-d(\tau_i))$ with $\tau_0 = 0$.

Taking into account that $w(\tau_i) = y(\tau_i) + v(\tau_i)$ for all $\tau_i \in \pi$ with $i \geq 1$ and that $\dot{w}(t) = L_f h(\Psi(z(t)))$, $t \in [\tau_i, \tau_{i+1})$, we get for all $t \in [\tau_i, \tau_{i+1}) \cap [0, t_{\max})$ with $i \geq 1$:

$$|w(t) - h(x(t))| = \left| y(\tau_i) + v(\tau_i) + \int_{\tau_i}^t L_f h(\Psi(z(s))) ds - y(t) \right| = \left| v(\tau_i) + \int_{\tau_i}^t L_f h(\Psi(z(s))) ds - \int_{\tau_i}^t L_f h(x(s)) ds \right|$$

The above equality in conjunction with the fact that $0 < \tau_{i+1} - \tau_i \leq r$ and estimate (3.5) implies for all $t \in [\tau_i, \tau_{i+1}) \cap [0, t_{\max})$ with $i \geq 1$:

$$\begin{aligned} |w(t) - h(x(t))| &\leq r \sup_{\tau_i \leq s \leq t} |L_f h(\Psi(z(s))) - L_f h(x(s))| + |v(\tau_i)| \leq r\bar{\sigma}(|x_0| + |z_0|, \tau_i) + rK \sup_{0 \leq \tau \leq t} |w(\tau) - h(x(\tau))| + |v(\tau_i)| \\ &\leq \sigma_1(|x_0| + |z_0|, t) + rK \sup_{0 \leq \tau \leq t} |w(\tau) - h(x(\tau))| + \sup_{0 \leq \tau \leq t} |v(\tau)| \end{aligned} \tag{3.7}$$

where $\sigma_1(s, t) := r\bar{\sigma}(s, t - r)$ for $t \geq r$ and $\sigma_1(s, t) := \exp(r - t) r \bar{\sigma}(s, 0)$ for $t < r$.

On the other hand taking into account that $w(\tau_0) = w(0) = w_0$ and that $\dot{w}(t) = L_f h(\Psi(z(t)))$, $t \in [0, \tau_1)$, we get for all $t \in [0, \tau_1) \cap [0, t_{\max})$:

$$|w(t) - h(x(t))| \leq |w_0 - h(x_0)| + \left| \int_0^t L_f h(\Psi(z(s))) ds - \int_0^t L_f h(x(s)) ds \right|$$

Continuity of $h$ in conjunction with the fact that $h(0) = 0$ implies the existence of a function $\rho \in K_\infty$ such that

$$|w - h(x)| \leq \rho(|w| + |x|), \quad \forall (x, w) \in \Re^n \times \Re \tag{3.8}$$

The previous inequalities in conjunction with the fact that $0 < \tau_{i+1} - \tau_i \leq r$ and estimate (3.5) imply for all $t \in [0, \tau_1) \cap [0, t_{\max})$:

$$\begin{aligned}
|w(t) - h(x(t))| &\leq \rho(|w_0| + |x_0|) + r \sup_{0 \leq s \leq t} |L_f h(\Psi(z(s))) - L_f h(x(s))| \\
&\leq \rho(|w_0| + |x_0|) + r\bar{\sigma}(|x_0| + |z_0|, 0) + rK \sup_{0 \leq \tau \leq t} |w(\tau) - h(x(\tau))| \\
&\leq \sigma_2(|x_0| + |z_0| + |w_0|, t) + rK \sup_{0 \leq \tau \leq t} |w(\tau) - h(x(\tau))|
\end{aligned} \tag{3.9}$$

where $\sigma_2(s, t) := r\bar{\sigma}(s, t - r) + \rho(s)\exp(r - t)$ for $t \geq r$ and $\sigma_2(s, t) := [r\bar{\sigma}(s, 0) + \rho(s)]\exp(r - t)$ for $t < r$. Notice that $\sigma_2 \in KL$. Combining estimates (3.7) and (3.9) we conclude that the following estimate holds for all $t \in [0, t_{\max})$:

$$|w(t) - h(x(t))| \leq \sigma_2(|x_0| + |z_0| + |w_0|, t) + rK \sup_{0 \leq \tau \leq t} |w(\tau) - h(x(\tau))| + \sup_{0 \leq \tau \leq t} |v(\tau)| \tag{3.10}$$

Using (3.8), (3.10) and (2.1) we obtain for all $t \in [0, t_{\max})$:

$$|w(t)| \leq \rho(\mu^2(t)) + \rho((a(|x_0|))^2) + \sigma_2(|x_0| + |z_0| + |w_0|, 0) + rK \sup_{0 \leq \tau \leq t} |w(\tau) - h(x(\tau))| + \sup_{0 \leq \tau \leq t} |v(\tau)|$$

The above inequality in conjunction with (3.6) gives for all $t \in [0, t_{\max})$:

$$|x(t)| + |z(t)| + |w(t)| \leq \phi(t) + \tilde{a}(|x_0| + |z_0| + |w_0|) + \sup_{0 \leq \tau \leq t} q(|w(\tau) - h(x(\tau))|) + \sup_{0 \leq \tau \leq t} |v(\tau)| \tag{3.11}$$

where $\phi(t) := \frac{1}{2}\mu^2(t) + \rho(\mu^2(t))$, $\tilde{a}(s) := \rho((a(s))^2) + (a(s))^2 + \sigma_2(s, 0)$ and $q(s) := (p(s))^2 + rKs$.

Using (3.10) and the fact that $rK < 1$, we obtain:

$$\sup_{0 \leq \tau < t_{\max}} |w(\tau) - h(x(\tau))| \leq \frac{1}{1 - rK} \sigma_2(|x_0| + |z_0| + |w_0|, 0) + \frac{1}{1 - rK} \sup_{0 \leq \tau \leq t_{\max}} |v(\tau)| \tag{3.12}$$

Exploiting (2.1), (3.6), (3.8) and (3.12) and the Boundedness-Implies-Continuation property for system (3.2) we may conclude that $t_{\max} = +\infty$. It follows that all the above inequalities hold for all $t \geq 0$. Moreover, taking into account that system (3.2) is autonomous, we may utilize (2.1), (3.6), (3.8) and (3.12) in order to show that $0 \in \Re^n \times \Re^k \times \Re$ is a robust equilibrium point from the input $v \in L_{loc}^\infty(\Re^+; \Re)$, i.e., for every $\varepsilon > 0$, $T \in \Re^+$ there exists $\delta := \delta(\varepsilon, T) > 0$ such that for all $(x_0, z_0, d, v) \in \Re^n \times \Re^k \times L_{loc}^\infty(\Re^+; \Re^+) \times L_{loc}^\infty(\Re^+; \Re)$, with $|x_0| + |z_0| + \sup_{t \geq 0}|v(t)| < \delta$ it holds that the solution $(x(t), z(t), w(t))$ of (3.2) with initial condition $(x(0), z(0), w(0)) = (x_0, z_0, w_0)$ corresponding to $(d, v) \in L_{loc}^\infty(\Re^+; \Re^+) \times L_{loc}^\infty(\Re^+; \Re)$ exists for all $t \in [t_0, t_0 + T]$ and

$$\sup\{|(x(t), z(t), w(t))|;\ t \in [0, T]\} < \varepsilon$$

Using (3.10), (3.11), Theorem 3.1 in [17] (Small-Gain Theorem for hybrid systems) in conjunction with Remarks 3.2 and 3.6 in [17], inequality $rK < 1$ and the fact that system (3.2) is autonomous, we conclude that system (3.2) is Robustly Forward Complete and there exist functions $\sigma_3 \in KL$, $\bar{\gamma} \in N$ such that for every $(x_0, z_0, w_0, d, v) \in \Re^n \times \Re^k \times \Re \times L_{loc}^\infty(\Re^+; \Re^+) \times L_{loc}^\infty(\Re^+; \Re)$ the solution $(x(t), z(t), w(t))$ of (3.2) with initial condition $(x(0), z(0), w(0)) = (x_0, z_0, w_0)$ corresponding to $d \in L_{loc}^\infty(\Re^+; \Re^+)$, $v \in L_{loc}^\infty(\Re^+; \Re)$ satisfies:

$$|w(t) - h(x(t))| \leq \sigma_3\left(|x_0| + |z_0| + |w_0|, t\right) + \sup_{0 \leq \tau \leq t} \bar{\gamma}(|v(\tau)|), \quad \forall t \geq 0 \qquad (3.13)$$

Using (3.4), (3.11), (3.13), Theorem 3.1 in [17] (Small-Gain Theorem for hybrid systems) in conjunction with Remark 3.2 in [17] and the fact that system (3.2) is autonomous, we conclude that there exist functions $\sigma \in KL$, $\tilde{\gamma} \in N$ such that for every $(x_0, z_0, w_0, d, v) \in \Re^n \times \Re^k \times \Re \times L^\infty_{loc}(\Re^+; \Re^+) \times L^\infty_{loc}(\Re^+; \Re)$ the solution $(x(t), z(t), w(t))$ of (3.2) with initial condition $(x(0), z(0), w(0)) = (x_0, z_0, w_0)$ corresponding to $d \in L^\infty_{loc}(\Re^+; \Re^+)$, $v \in L^\infty_{loc}(\Re^+; \Re)$ satisfies (3.3). The proof is complete. ◁

**Remark 3.3:** It should be noted that inequality (3.1) is a very conservative condition, which is rarely satisfied for general nonlinear observers. The following section provides classes of nonlinear observers that satisfy condition (3.1). On the other hand, the reader should notice that condition (3.1) may be satisfied by simply redefining the output map of system (1.1). For example, if a function $b \in C^2(\Re; \Re)$ with $b(0) = 0$, $b(\Re) = \Re$ and $\frac{db}{dx}(x) \neq 0$ for all $x \in \Re$ can be found so that the following system:

$$\dot{z} = \tilde{F}(z, \tilde{y}), z \in \Re^k$$
$$\hat{x} = \Psi(z), \hat{x} \in \Re^n \qquad (3.14)$$

where $\tilde{F}(z, \tilde{y}) = F(z, b^{-1}(\tilde{y}))$, $b^{-1} \in C^1(\Re; \Re)$ is the inverse function of $b \in C^2(\Re; \Re)$, is a robust observer with respect to measurement errors for system (1.1) with output $\tilde{y} = b(h(x))$ instead of $y = h(x)$ and there exists a constant $K \geq 0$ and a function $\bar{\sigma} \in KL$ such that for every $(x_0, z_0) \in \Re^n \times \Re^k$ and $e \in L^\infty_{loc}(\Re^+; \Re)$, the solution $(x(t), z(t))$ of

$$\dot{x} = f(x)$$
$$\dot{z} = \tilde{F}(z, b(h(x)) + e) \qquad (3.15)$$
$$\hat{x} = \Psi(z)$$

satisfies the following estimate:

$$\left|\frac{db}{ds}(h(\Psi(z(t))))L_f h(\Psi(z(t))) - \frac{db}{ds}(h(x(t)))L_f h(x(t))\right| \leq \bar{\sigma}(|x_0| + |z_0|, t) + K \sup_{0 \leq \tau \leq t}|e(\tau)|, \quad \forall t \geq 0 \quad (3.16)$$

then a sampled-data observer for system (1.1) with output $\tilde{y} = b(h(x))$ instead of $y = h(x)$ exists and Theorem 3.1 applies.

## 4. Applications

In this section we present the application of Theorem 3.1 to two classes of systems: (i) triangular globally Lipschitz systems and (ii) linear detectable systems, leading to concrete sampled-data observer designs.

a) <u>Triangular Globally Lipschitz Systems</u>:

Consider the system

$$\dot{x}_1 = f_1(x_1) + x_2$$
$$\vdots$$
$$\dot{x}_{n-1} = f_{n-1}(x_1, \ldots, x_{n-1}) + x_n \qquad (4.1)$$
$$\dot{x}_n = f_n(x_1, \ldots, x_n)$$
$$y = x_1$$

where $f_i : \Re^i \to \Re$ ($i = 1, \ldots, n$) with $f_i(0) = 0$ ($i = 1, \ldots, n$) are globally Lipschitz functions, i.e., there exists a constant $L \geq 0$ such that the following inequalities hold for $i = 1, \ldots, n$:

$$|f_i(x_1,...,x_i) - f_i(z_1,...,z_i)| \leq L|(x_1-z_1,...,x_i-z_i)|, \quad \forall (x_1,...,x_i) \in \Re^i, \quad \forall (z_1,...,z_i) \in \Re^i \quad (4.2)$$

The reader should notice that all linear observable systems can be written in the form (4.1) with $f_i : \Re^i \to \Re$ ($i = 1,...,n$) being linear functions. Notice that systems of the form (4.1) are Robustly Forward Complete and satisfy hypothesis (H), since for every $x_0 \in \Re^n$ the solution of (4.1) with initial condition $x(0) = x_0$ satisfies the estimate:

$$|x(t)| \leq \exp(ct)|x_0|, \quad \forall t \geq 0 \quad (4.3)$$

where $c := nL + n - 1$. Inequality (4.3) is obtained by evaluating the derivative of the function $W(x) = \frac{1}{2}|x|^2$ along the solutions of (4.1) and using inequalities (4.2).

A high-gain observer design is described in [12]: first a vector $k = (k_1,...,k_n)' \in \Re^n$ is found so that the matrix $(A + kc') \in \Re^{n \times n}$ is Hurwitz, where $c := (1,0,...,0)' \in \Re^n$ and $A \in \Re^{n \times n}$ is the matrix $A = \{a_{i,j} : i = 1,...,n, j = 1,...,n\}$ with $a_{i,i+1} = 1$ for $i = 1,...,n-1$ and $a_{i,j} = 0$ if otherwise. The existence of the required vector $k = (k_1,...,k_n)' \in \Re^n$ is guaranteed since the pair of matrices $(A,c)$ is observable. The proposed observer is of the form:

$$\begin{aligned}
\dot{z}_i &= f_i(z_1,...,z_i) + z_{i+1} + \theta^i k_i(c'z - y), \quad i = 1,...,n-1 \\
\dot{z}_n &= f_n(z_1,...,z_n) + \theta^n k_n(c'z - y) \\
\hat{x} &= z \\
z &= (z_1,...,z_n)' \in \Re^n
\end{aligned} \quad (4.4)$$

where $\theta \geq 1$ is a constant sufficiently large. The proof is based on the quadratic error Lyapunov function $V(e) := e'\Delta_\theta^{-1} P \Delta_\theta^{-1} e$, where $e := z - x$, $\Delta_\theta := diag(\theta,\theta^2,...,\theta^n)$ and $P \in \Re^{n \times n}$ is a symmetric positive definite matrix that satisfies $P(A+kc') + (A'+ck')P + 2\mu I \leq 0$ for certain constant $\mu > 0$ (see [12] for details). Specifically, using the identities $\Delta_\theta^{-1} A = \theta A \Delta_\theta^{-1}$, $c' = \theta c' \Delta_\theta^{-1}$ and the inequalities $\theta^{-i}|f_i(x_1+e_1,...,x_i+e_i) - f_i(x_1,...,x_i)| \leq L|\Delta_\theta^{-1} e|$ for $i = 1,...,n$ and all $(x,e) \in \Re^n \times \Re^n$ (which follow from (4.2)), we get for $\theta \geq \max\left(1, \frac{2|P|L\sqrt{n}}{\mu}\right)$ and for all $(x,e) \in \Re^n \times \Re^n$:

$$\dot{e} = (A + \Delta_\theta kc')e + g(x,e) \quad (4.5)$$

$$\begin{aligned}
\nabla V(e)[(A + \Delta_\theta kc')e + g(x,e)] &\leq -2\theta\mu|\Delta_\theta^{-1} e|^2 + 2|\Delta_\theta^{-1} e||P||\Delta_\theta^{-1} g(x,e)| \leq \\
-2\theta\mu|\Delta_\theta^{-1} e|^2 + 2|\Delta_\theta^{-1} e|^2|P|L\sqrt{n} &\leq -\theta\mu|\Delta_\theta^{-1} e|^2 \leq -\frac{\theta\mu}{|P|}V(e)
\end{aligned} \quad (4.6)$$

where $g(x,e) = (f_1(x_1+e_1) - f_1(x_1),...,f_n(x+e) - f_n(x))'$. Using the inequality $2e'\Delta_\theta^{-1} Pkv \leq \frac{\theta\mu}{2|P|}V(e) + \frac{2|P|^2|k|^2}{\theta\mu}|v|^2$ (which holds for all $(e,v) \in \Re^n \times \Re$) we obtain the following inequality for all $(x,e,v) \in \Re^n \times \Re^n \times \Re$:

$$\nabla V(e)[(A + \Delta_\theta kc')e + \Delta_\theta kv + g(x,e)] \leq -\frac{\theta\mu}{2|P|}V(e) + \frac{2|P|^2|k|^2}{\theta\mu}|v|^2 \quad (4.7)$$

Inequality (4.7) implies that for all $(x_0, z_0, v) \in \Re^n \times \Re^n \times L_{loc}^\infty(\Re^+; \Re)$ the solution of (4.1) with

$$\dot{z}_i = f_i(z_1,...,z_i) + z_{i+1} + \theta^i k_i (c'z - c'x + v), \ i = 1,...,n-1$$
$$\dot{z}_n = f_n(z_1,...,z_n) + \theta^n k_n (c'z - c'x + v)$$
$$\hat{x} = z \quad (4.8)$$
$$z = (z_1,...,z_n)' \in \Re^n$$

and initial condition $(x(0), z(0)) = (x_0, z_0) \in \Re^n \times \Re^n$ corresponding to $v \in L_{loc}^\infty (\Re^+; \Re)$ satisfies the estimates:

$$|\hat{x}(t) - x(t)| \leq \theta^{n-1} \sqrt{\frac{K_2}{K_1}} \exp\left(-\frac{\theta \mu}{4|P|}\right) |z_0 - x_0| + 2\theta^{n-1} \frac{|P\|k|}{\mu} \sqrt{\frac{K_2}{K_1}} \sup_{0 \leq \tau \leq t} |v(\tau)|, \ \forall t \geq 0 \quad (4.9)$$

$$|z_i(t) - x_i(t)| \leq \theta^{i-1} \sqrt{\frac{K_2}{K_1}} \exp\left(-\frac{\theta \mu}{4|P|}\right) |z_0 - x_0| + 2\theta^{i-1} \frac{|P\|k|}{\mu} \sqrt{\frac{K_2}{K_1}} \sup_{0 \leq \tau \leq t} |v(\tau)|, \ \forall t \geq 0, \ i = 1,...,n \quad (4.10)$$

where $K_1, K_2 > 0$ are constants such that $K_1 |x|^2 \leq x'Px \leq K_2 |x|^2$ for all $x \in \Re^n$. It follows from (4.9) and (4.3) that system (4.4) is a robust observer for system (4.1) with respect to measurement errors. Moreover, using inequality (4.2) for $i = 1$ and (4.10) we obtain that for all $(x_0, z_0, v) \in \Re^n \times \Re^n \times L_{loc}^\infty (\Re^+; \Re)$ the solution of (4.1) with (4.8) and initial condition $(x(0), z(0)) = (x_0, z_0) \in \Re^n \times \Re^n$ corresponding to $v \in L_{loc}^\infty (\Re^+; \Re)$ satisfies the estimate:

$$|f_1(z_1(t)) + z_2(t) - f_1(x_1(t)) - x_2(t)| \leq (L+\theta) \sqrt{\frac{K_2}{K_1}} \exp\left(-\frac{\theta \mu}{4|P|}\right) |z_0 - x_0| + 2(L+\theta) \frac{|P\|k|}{\mu} \sqrt{\frac{K_2}{K_1}} \sup_{0 \leq \tau \leq t} |v(\tau)|, \ \forall t \geq 0 \quad (4.10)$$

It follows that the following system

$$\dot{z}_i(t) = f_i(z_1(t),...,z_i(t)) + z_{i+1}(t) + \theta^i k_i (c'z(t) - w(t)), \ i = 1,...,n-1$$
$$\dot{z}_n(t) = f_n(z_1(t),...,z_n(t)) + \theta^n k_n (c'z(t) - w(t))$$
$$\dot{w}(t) = f_1(z_1(t)) + z_2(t) \quad , \quad t \in [\tau_i, \tau_{i+1})$$
$$w(\tau_{i+1}) = y(\tau_{i+1}) \quad (4.11)$$
$$\hat{x} = z$$
$$z = (z_1,...,z_n)' \in \Re^n$$

is a robust sampled-data observer for system (4.1) with respect to measurement errors provided that the upper diameter of the sampling partition $r > 0$ satisfies the inequality:

$$2r(L+\theta) \frac{|P\|k|}{\mu} \sqrt{\frac{K_2}{K_1}} < 1 \quad (4.12)$$

Notice that since $\theta \geq \max\left(1, \frac{2|P|L\sqrt{n}}{\mu}\right)$, it follows from (4.12) that the upper diameter of the sampling partition $r > 0$ must necessarily be less than $\frac{\mu^2}{2|P\|k|(L\mu + \max(\mu, 2|P|L\sqrt{n}))} \sqrt{\frac{K_1}{K_2}}$.

b) <u>Linear detectable systems</u>:

Consider
$$\dot{x} = Ax$$
$$y = c'x \quad (4.13)$$
$$x \in \Re^n, y \in \Re$$

there exists a vector $k \in \Re^n$ such that the matrix $(A+kc')$ is Hurwitz. Consequently, there exists a positive definite symmetric matrix $P \in \Re^{n \times n}$ such that the symmetric matrix $P(A+kc')+(A'+ck')P$ is negative definite and there exist constants $\mu, \gamma > 0$ such that

$$x'P(A+kc')x + x'(A'+ck')Px + 2x'Pkv \leq -2\mu x'Px + \gamma |v|^2, \quad \forall (x,v) \in \Re^n \times \Re \qquad (4.14)$$

It follows that for every $(x_0, z_0, v) \in \Re^n \times \Re^n \times L^\infty_{loc}(\Re^+; \Re)$ the solution of (4.13) with

$$\dot{z} = Az + k(c'z - y + v)$$
$$\hat{x} = z \qquad (4.15)$$

and initial condition $(x(0), z(0)) = (x_0, z_0) \in \Re^n \times \Re^n$ corresponding to $v \in L^\infty_{loc}(\Re^+; \Re)$ satisfies the estimates:

$$|z(t) - x(t)| \leq \exp(-\mu t)\sqrt{\frac{K_2}{K_1}}|z_0 - x_0| + \sqrt{\frac{\gamma}{2\mu K_1}} \sup_{0 \leq \tau \leq t} |v(\tau)|, \quad \forall t \geq 0 \qquad (4.16)$$

$$|c'Az(t) - c'Ax(t)| \leq |c'A|\exp(-\mu t)\sqrt{\frac{K_2}{K_1}}|z_0 - x_0| + |c'A|\sqrt{\frac{\gamma}{2\mu K_1}} \sup_{0 \leq \tau \leq t} |v(\tau)|, \quad \forall t \geq 0 \qquad (4.17)$$

where $K_1, K_2 > 0$ are constants such that $K_1|x|^2 \leq x'Px \leq K_2|x|^2$ for all $x \in \Re^n$. It follows from Theorem 3.1 that the following system

$$\dot{z}(t) = Az(t) + k(c'z(t) - w(t))$$
$$\dot{w}(t) = c'Az(t), \quad t \in [\tau_i, \tau_{i+1})$$
$$w(\tau_{i+1}) = y(\tau_{i+1})$$
$$\hat{x} = z \qquad (4.18)$$

is a robust sampled-data observer for system (4.13) with respect to measurement errors provided that the upper diameter of the sampling partition $r > 0$ satisfies the inequality:

$$r|c'A|\sqrt{\frac{\gamma}{2\mu K_1}} < 1 \qquad (4.19)$$

As an example, consider the linear oscillator:

$$\dot{x}_1 = x_2$$
$$\dot{x}_2 = -4x_1 \qquad (4.20)$$
$$y = x_1$$

Inequality (4.14) holds with $P = \frac{1}{2}\begin{bmatrix} 5 & -2 \\ -2 & 1 \end{bmatrix}$, $k = (-4, 0)'$, $c = (1, 0)'$, $\mu = 1$, $\gamma = \frac{64}{3}$. Using the fact that $(3 - 2\sqrt{2})|x|^2 \leq 2x'Px \leq 7|x|^2$, for all $x \in \Re^n$, we conclude from (4.19) that the following system

$$\dot{z}_1 = z_2 - 4(z_1 - w)$$
$$\dot{z}_2 = -4z_1$$
$$\dot{w} = z_2, \quad t \in [\tau_i, \tau_{i+1}) \qquad (4.21)$$
$$w(\tau_{i+1}) = y(\tau_i)$$
$$\hat{x} = (z_1, z_2)$$

is a sampled-data observer for (4.20) provided that $r < 0.089$. Numerical results are presented for sampling partition $\pi = \{\tau_i\}_{i=0}^{\infty}$ satisfies $\tau_i = ir$, $i = 0,1,2,...$ with $r = 0.081$. For comparison purpose we compare the performance of the sampled-data observer (4.21) with the continuous-time observer:

$$\begin{aligned}\dot{z}_1 &= z_2 - 4(z_1 - y) \\ \dot{z}_2 &= -4z_1 \\ \hat{x} &= (z_1, z_2)\end{aligned} \quad (4.22)$$

as well as with the discrete implementation of the continuous-time observer (4.22), where the output $y(t)$ for $t \in [\tau_i, \tau_{i+1})$ is simply replaced by the most recent measurement $y(\tau_i)$:

$$\begin{aligned}\dot{z}_1 &= z_2 - 4(z_1 - y(\tau_i)), \, t \in [\tau_i, \tau_{i+1}) \\ \dot{z}_2 &= -4z_1 \\ \hat{x} &= (z_1, z_2)\end{aligned} \quad (4.23)$$

Notice that system (4.23) is the usual implementation of the continuous-time observer (4.22). Figures 2-4 present the evolution of the error variable $e_2 = x_2 - z_2$ with initial conditions $x_1(0) = 0$, $x_2(0) = 2$, $z_1(0) = z_2(0) = 1$ and $w(0) = 0$. As expected the error variable $e_2 = x_2 - z_2$ does not converge to zero for (4.23): it presents a persistent oscillation with approximate amplitude 0.15.

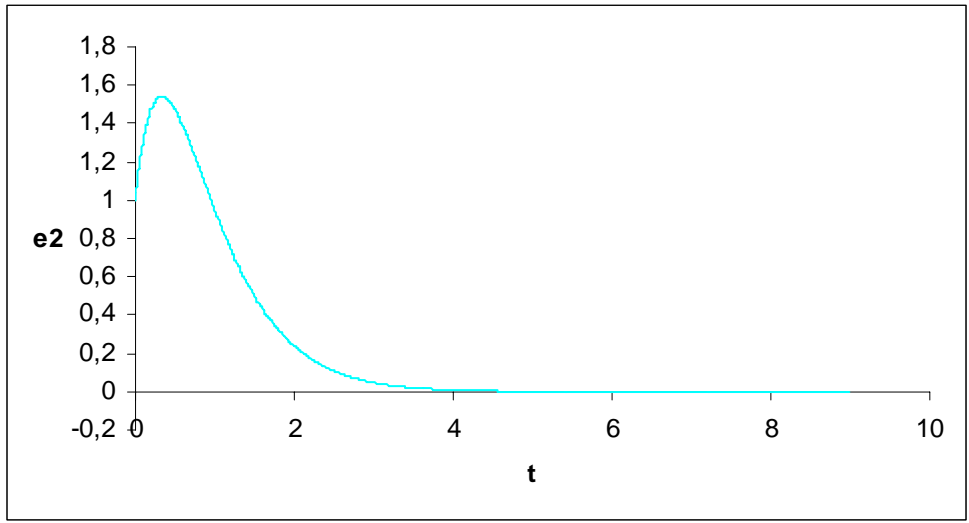

**Figure 2:** Evolution of $e_2 = x_2 - z_2$ for the continuous-time observer (4.22)

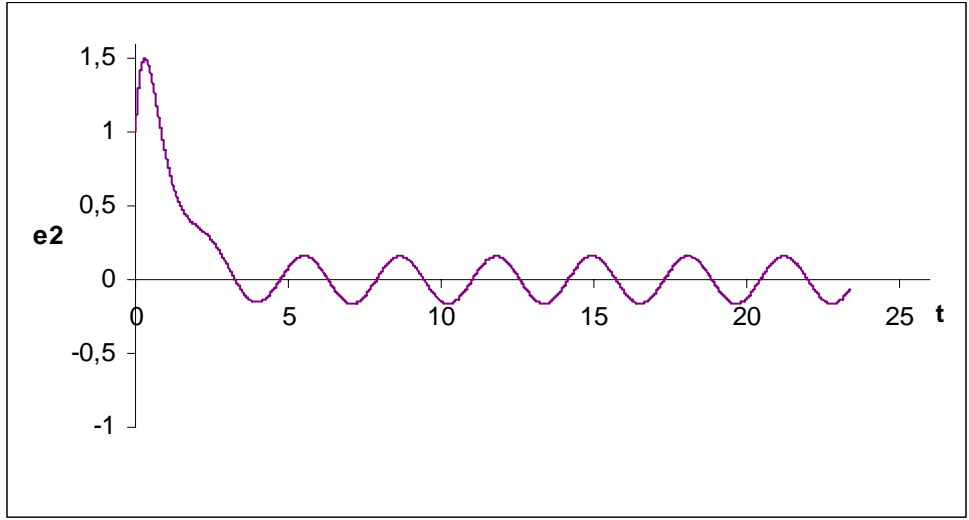

**Figure 3:** Evolution of $e_2 = x_2 - z_2$ for (4.23), $r = 0.081$

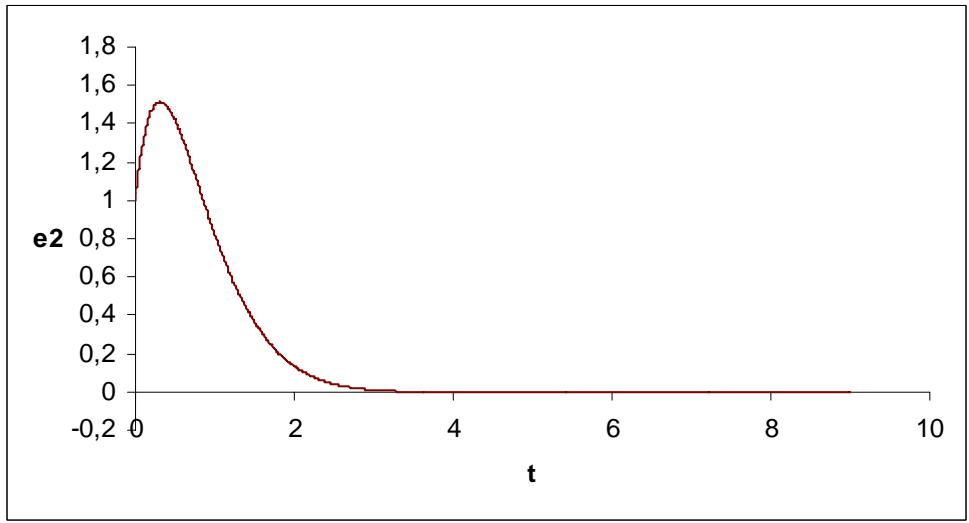

**Figure 4:** Evolution of $e_2 = x_2 - z_2$ for the sampled-data observer (4.21), $r = 0.081$

It should be noted here that condition (4.19) is very conservative. Indeed, simulations show that system (4.21) is a sampled-data observer for uniform sampling partitions $\pi = \{\tau_i\}_{i=0}^{\infty}$, $\tau_i = ir$, $i = 0,1,2,...$ with $r > 0.081$. Figures 5-6 present the evolution of the error variable $e_2 = x_2 - z_2$ with initial conditions $x_1(0) = 0$, $x_2(0) = 2$, $z_1(0) = z_2(0) = 1$ and $w(0) = 0$ for uniform sampling partition $\pi = \{\tau_i\}_{i=0}^{\infty}$, $\tau_i = ir$, $i = 0,1,2,...$ with $r = 0.45$. It is clear that system (4.23) gives a very disappointing state estimation, while the sampled-data observer provides a reliable estimation of the state.

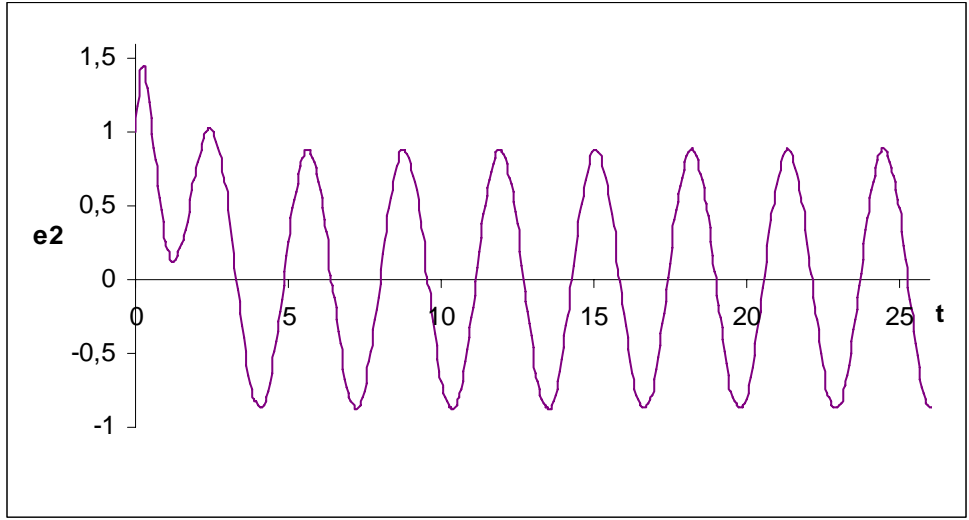

**Figure 5:** Evolution of $e_2 = x_2 - z_2$ for (4.23), $r = 0.45$

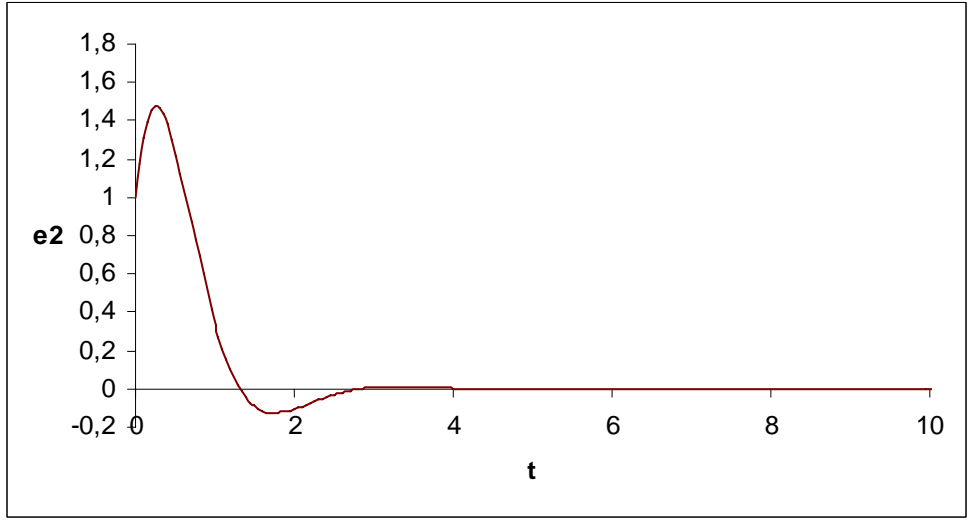

**Figure 6:** Evolution of $e_2 = x_2 - z_2$ for the sampled-data observer (4.21), $r = 0.45$

Having compared the proposed observer (4.21) with the usual implementation of the continuous-time observer (4.22), we now compare the proposed observer (4.21) with discrete-time observer designs. Assuming uniform sampling partition $\pi = \{\tau_i\}_{i=0}^{\infty}$ with $\tau_i = iT$, $i = 0,1,2,...$, where $T > 0$ denotes the sampling period, we obtain the following discrete-time model for (4.20):

$$x_{k+1} = A(T)x_k = \begin{bmatrix} \cos(2T) & \frac{1}{2}\sin(2T) \\ -2\sin(2T) & \cos(2T) \end{bmatrix} x_k \qquad (4.24)$$
$$y_k = c'x_k = \begin{bmatrix} 1 & 0 \end{bmatrix} x_k$$

where $x_k = x(kT) \in \Re^2$. A discrete-time design is the selection of a vector $L = (L_1, L_2)' \in \Re^2$ so that the eigenvalues of the matrix $A(T) + Lc'$ are placed strictly inside the unit ball on the complex plane. As long as $\frac{2T}{\pi}$ is not an integer, this task can be accomplished and we can place the eigenvalues of the matrix $A(T) + Lc'$ at any desired positions. The observer is implemented by:

$$z_{k+1} = A(T)z_k + L(c'z_k - y_k) = \begin{bmatrix} \cos(2T) + L_1 & \frac{1}{2}\sin(2T) \\ -2\sin(2T) + L_2 & \cos(2T) \end{bmatrix} z_k - \begin{bmatrix} L_1 \\ L_2 \end{bmatrix} y_k \quad (4.25)$$

Of course, an additional rule must be given for the estimation of the state variables between two consecutive sampling instances, but in what follows intersampling behavior will be ignored. Suppose now that the sampling partition $\pi = \{\tau_i\}_{i=0}^{\infty}$ does not necessarily satisfy $\tau_i = iT$, $i = 0,1,2,...$. If the observer design is based on the hypothesis of constant sampling period $T > 0$ then the error $e_k = z_k - x(\tau_k)$ satisfies the equation:

$$e_{k+1} = (A(T) + Lc')e_k + (A(T) - A(\tau_{k+1} - \tau_k))x(\tau_k) \quad (4.26)$$

where the equations $x(\tau_{k+1}) = A(\tau_{k+1} - \tau_k)x(\tau_k)$ and $y_k = c'x(\tau_k)$ have been used for the derivation of the above equation. The above equation in conjunction with the fact that the matrix $A(T) + Lc'$ is a Schur matrix implies the existence of constants $K, \gamma, c > 0$ such that the error satisfies the following estimate for all $k \geq 1$:

$$|e_k| \leq K \exp(-ck)|e_0| + \gamma \max_{0 \leq j < k} |(A(T) - A(\tau_{j+1} - \tau_j))x(\tau_j)| \quad (4.27)$$

It should be clear from the above estimates that we cannot conclude that $e_k \to 0$. Indeed, Figure 7 shows the long-time behavior of the error variable $e_2(\tau_k) = \begin{bmatrix} 0 & 1 \end{bmatrix}(z_k - x(\tau_k))$ for sampling partition $\pi = \{\tau_i\}_{i=0}^{\infty}$ which satisfies $\tau_i = ir$, $i = 0,1,2,...$ with $r = 0.081$, assumed sampling period $T = 0.075$ (i.e., 8% error in sampling period) and initial condition $x_1(0) = 0$, $x_2(0) = 2$, $z_0 = (1,1)' \in \Re^2$. The observer gains $L = (L_1, L_2)' \in \Re^2$ have been selected so that both eigenvalues of the matrix $A(T) + Lc'$ are equal to 0.8 ($L_1 = -0.37754$ and $L_2 = -0.17804$). It is clear that the error does not converge to zero but instead presents an oscillation of amplitude greater than 0.21. The reader should compare Figure 7 with Figure 4. Similar results are obtained for different selections of the observer gains $L = (L_1, L_2)' \in \Re^2$; indeed, in all cases the error does not converge to zero but instead presents an oscillation of amplitude greater than 0.15.

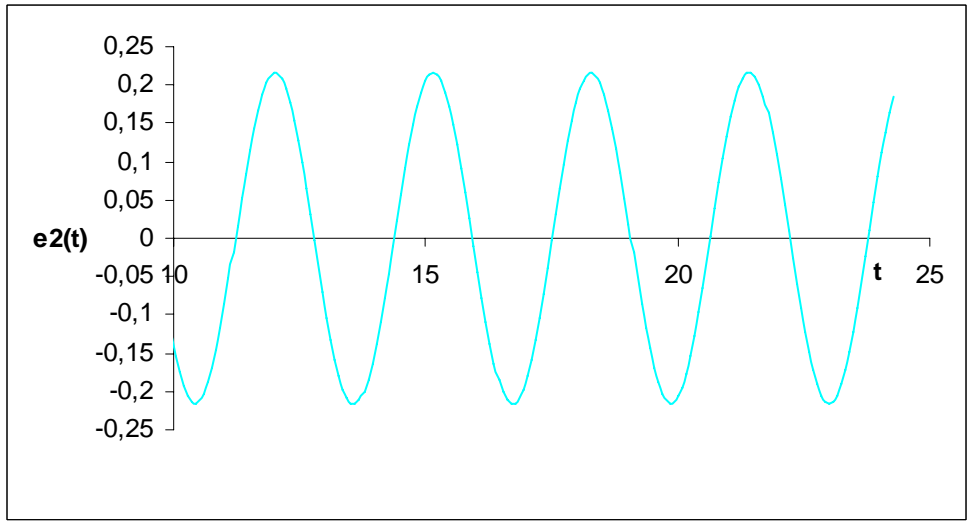

**Figure 7:** Long-time behavior of the error for the discrete-time observer (4.25)

As remarked in the Introduction, discrete-time observer design seems to be sensitive to perturbations of the sampling schedule; on the other hand the proposed observer (4.21) is insensitive to perturbations of the sampling schedule.

# 5. Concluding Remarks

The present work developed a design method for nonlinear sampled-data observers based on an available continuous-time design, coupled with an inter-sample output predictor.

In addition to being intuitively meaningful, key attractive features of the proposed sampled-data observer include that
  1) it provides easily checkable sufficient conditions for robustness with respect to measurement errors.
  2) it provides an explicit formula for estimating the maximum allowable sampling period.
  3) it provides explicit bounds for the estimation error between sampling instants.
  4) it provides robustness with respect to perturbations of the sampling schedule.
The theory leads to particularly simple results for the linear case and for the globally Lipschitz case.

At present, the theoretical formulation strictly refers to global observers. Future work will study a local formulation of the theory, to enlarge the range of applicability of the sampled-data observer.